\documentclass[english]{article}
\usepackage[T2A]{fontenc}
\usepackage[utf8]{inputenc}
\usepackage[a4paper]{geometry}
\usepackage{babel}
\usepackage{amsthm}
\usepackage{amsmath}
\usepackage{amssymb}
\usepackage[unicode=true]
 {hyperref}

\makeatletter
\theoremstyle{plain}
\newtheorem{thm}{\protect\theoremname}
\theoremstyle{definition}
\newtheorem{defn}[thm]{\protect\definitionname}
\theoremstyle{plain}
\newtheorem{lem}[thm]{\protect\lemmaname}
\newenvironment{lyxlist}[1]
{\begin{list}{}
{\settowidth{\labelwidth}{#1}
 \setlength{\leftmargin}{\labelwidth}
 \addtolength{\leftmargin}{\labelsep}
 }}
{\end{list}}

\newcommand{\vestnikonly}[1]{}
\newcommand{\novestnikonly}[1]{#1}
\vestnikonly{
\usepackage{allerree}

\usepackage{nameref}
\usepackage{url}
\let\href\undefined
\newcommand{\href}[2]{\url{#1}}

\newcounter{countergvtheorems}

\renewenvironment{lem}
{\refstepcounter{countergvtheorems}\begin{theorem}{Лемма  \arabic{countergvtheorems}.}}
{\end{theorem}\par{\bf Доказательство.}}

\renewenvironment{thm}
{\refstepcounter{countergvtheorems}\begin{theorem}{Теорема \arabic{countergvtheorems}.}}
{\end{theorem}\par{\bf Доказательство.}}

\renewenvironment{defn}
{\par \refstepcounter{countergvtheorems} {\bf Определение \arabic{countergvtheorems}.}}
{\par}
}
\novestnikonly{
\newcommand{\iabstract}[4]{
\begin{abstract}
#3
\footnote{\emph{Keywords:} #4}
\end{abstract}
}
}

\DeclareMathOperator{\HC}{HC}
\DeclareMathOperator{\HM}{HM}
\DeclareMathOperator{\LCM}{LCM}
\DeclareMathOperator{\poly}{poly}
\DeclareMathOperator{\SIG}{\mathcal{S}}

\AtBeginDocument{
  
}

\makeatother

\providecommand{\definitionname}{\inputencoding{latin9}Definition}
\providecommand{\lemmaname}{\inputencoding{latin9}Lemma}
\providecommand{\theoremname}{\inputencoding{latin9}Theorem}

\begin{document}
\global\long\def\GVWl{<_{\text{H}}}

\global\long\def\GVWg{>_{\text{\textnormal{H}}}}

\global\long\def\eqdef{\overset{\mathrm{_{def}}}{=}}

\global\long\def\equivdef{\overset{\mathrm{_{def}}}{\Leftrightarrow}}

\novestnikonly{

\title{Simple signature-based Groebner basis algorithm}

\author{Galkin Vasily\\
Moscow State University\\
email:\href{mailto: galkin-vv@yandex.ru}{ galkin-vv@yandex.ru}}

\maketitle
\providecommand{\definitionname}{Definition}
\providecommand{\lemmaname}{Lemma}
\providecommand{\theoremname}{Theorem} 
}
\vestnikonly{
\cleanbegin 
\def\udk{512}
\ltitle{ПРОСТОЙ ИТЕРАТИВНЫЙ АЛГОРИТМ ВЫЧИСЛЕНИЯ БАЗИСОВ ГРЁБНЕРА, ОСНОВАННЫЙ НА СИГНАТУРАХ} {В.\,В.~Галкин\footnote[1]{{\it Галкин Василий Витальевич} --- асп. каф. алгебры мех.-мат. ф-та МГУ, e-mail: galkin-vv@yandex.ru.}}
}
\iabstract
{Данная работа описывает алгоритм вычисления базисов Грёбнера, основанный на использовании отмеченных многочленов и идеях из алгоритма F5. Отличительными особенностями алгоритма по сравнению с аналогами являются простота самого алгоритма и доказательства его корректности, достигнутые без потери эффективности. Это позволяет создать простую реализацию, не уступающую более сложным аналогам по производительности}
{базис Грёбнера, алгоритм F5, отмеченные многочлены}
{This paper presents an algorithm for computing Groebner bases based upon labeled polynomials and ideas from the algorithm F5. The main highlights of this algorithm compared with analogues are simplicity both of the algorithm and of the its correctness proof achieved without loss of the efficiency. This leads to simple implementation which performance is in par with more complex analogues} {Groebner basis, F5 algorithm, labeled polynomials}

Consider polynomial ring $P=k[x_{1},\dots,x_{n}]$ over field $k$.
Also assume that monoid of its monomials $\mathbb{T}$ has a monomial
order$\prec$. A problem asking for a Gröbner basis can be stated
for any ideal$\left(f_{1},\dots,f_{l}\right)$ in this ring. One of
the approaches to the problem is using iterative method which computes
every step a basis for ideal $\left(f_{1},\dots,f_{i}\right),i=2\dots l$
based on the already computed for $\left(f_{1},\dots,f_{i-1}\right)$
basis $R_{i-1}$ and polynomial $f_{i}$. The algorithm described
in this paper is designed to perform one step of such computation.
So, the algorithm's input data consist of a some polynomial $f$ and
a polynomial set referred as $\left\{ g_{1},\dots,g_{m}\right\} $
which is Gröbner basis of ideal $I_{0}=\left(g_{1},\dots,g_{m}\right)$.
After finishing the algorithm should give the resulting polynomial
set $R$ being a Gröbner basis of ideal $I=\left(g_{1},\dots,g_{m},f\right)$.
The special cases $f=0\Rightarrow I=I_{0}$ and $\exists i\, g_{i}\in k\Rightarrow I=P$
are not interesting from the computational point of view, so the further
chapters assume that $f\neq0,\forall i\, g_{i}\notin k.$ The homogeneity
of input polynomials is not required unlike the F5 algorithm described
in \cite{FaugereF5}.

\section{Definitions}

Consider the set$\mathbb{T}_{0}=\mathbb{T}\cup\left\{ 0\right\} $
-- the monomial monoid extended by zero. The order $\prec$ can be
extended to $\mathbb{T}_{0}$ as $\prec_{0}$ with definition $\forall t\in\mathbb{T}\, t\succ_{0}0$
which keeps the well-orderness property. The notion of division also
can be extended to $\mathbb{T}_{0}$: $t_{1}|t_{2}\eqdef\exists t_{3}\, t_{1}t_{3}=t_{2}$.
For polynomial $p\in P,p\neq0$ the highest by $\prec$ monom and
coefficient are written as $\HM(p)\in\mathbb{T}$ and $\HC(p)\in k$.
For zero we define: $\HM(0)\eqdef0\in\mathbb{T}_{0}$, $\HC(0)\eqdef0\in k$.
The least common multiple of $t_{1},t_{2}\in\mathbb{T}$ is written
as $\LCM(t_{1},t_{2})\in\mathbb{T}$. In the following all definitions
are given for fixed $I_{0}$ and $f$:
\begin{defn}
The \emph{labeled polynomial} is a pair $h=(\sigma,p)\in\mathbb{T}_{0}\times P$,
that satisfies the correctness property: $\exists u\in P\,\HM(u)=\sigma,uf\equiv p\pmod{I_{0}}.$
Some terminology is extended to labeled polynomials. The highest monomial
is $\HM(h)\eqdef\HM(p)$ and coefficient is $\HC(h)\eqdef\HC(p)$.
Additionally the \emph{signature} is defined $\SIG(h)\eqdef\sigma$
and a notation is introduced for the polynomial -- second element
of pair: $\poly(h)\eqdef p$. The set of all labeled polynomials is
written as $H\subset\mathbb{T}_{0}\times P$. The trivial examples
of labeled polynomials are $\left(1,f\right)$ and $\left(0,g\right)$
for $g\in I_{0}$. Another labeled polynomial example is $\left(\HM(g),0\right)$
for $g\in I_{0}$. It satisfies correctness property because we can
take $u$ equal to $g$.\end{defn}
\begin{lem}
The product of $h\in H,t\in\mathbb{T}$ defined as $th\eqdef(t\sigma,tp)\in H,$
is correct.
\end{lem}
The correctness property is checked by directly finding $u$ for $th$.
\begin{defn}
If\emph{ $h'_{1},h_{2}\in H,t\in\mathbb{T}$} satisfy $\SIG(h'_{1})\succ_{0}\SIG(th_{2}),\HM(h'{}_{1})=\HM(th{}_{2})\neq0,$
then exists a \emph{signature-safe reduction $h'_{1}$ by $h_{2}$},
resulting in labeled polynomial $h_{1}\in H$, equal to: 
\[
h_{1}=\left(\SIG(h'_{1}),\poly(h'_{1})+Kt\poly(h_{2})\right),
\]
where the $K\in k$ is selected in a way to perform cancellation of
high coefficients, so we have $\HM(h{}_{1})\prec_{0}\HM(h'{}_{1})$.
Such reduction is equivalent to plain reduction with high term cancellation
extended with requirement for reductor's signature being smaller than
the signature of labeled polynomial being reduced. Like in previous
case the correctness check is performed directly.
\end{defn}
Let's introduce a partial order $\GVWl$ on $H$: 
\[
h_{1}=(\sigma_{1},p_{1})\GVWl h_{2}=(\sigma_{2},p_{2})\equivdef\HM(p_{1})\sigma_{2}\prec_{0}\HM(p_{2})\sigma_{1}.
\]

The elements with zero signature or zero high monomial are extremums:
\[
\forall\sigma_{1},\sigma_{2},p_{1},p_{2}\,\left(0,p_{1}\right)\not\GVWl\left(\sigma_{2},p_{2}\right),\,\left(\sigma_{1},0\right)\not\GVWg\left(\sigma_{2},p_{2}\right).
\]

\begin{lem}
Let $h_{1},h_{2}\in H,t\in\mathbb{T}$. Then $h_{1}\GVWg h_{2}\Leftrightarrow h_{1}\GVWg th_{2}$.
\end{lem}
Deduced from the fact that multiplying one of the compared labeled
polynomials by $t$ leads to multiplying by $t$ both sides in the
definition of$\GVWg$.
\begin{lem}
Let $h_{1},h_{2}\in H,\HM(h_{1})|\HM(h_{2}),\HM(h_{2})\neq0$. Then
signature-safe reduction $h_{2}$ by $h_{1}$ is possible iff$h_{1}\GVWg h_{2}$.
\end{lem}
Deduced from the fact that claims of both sides are equivalent to
$\SIG(h_{2})\succ_{0}\SIG(h_{1})\frac{\HM(h_{2})}{\HM(h_{1})}$.
\begin{lem}
Let $h_{1}\in H$ be a result of signature-safe reduction of $h_{1}'$
by some other polynomial. Then $h_{1}\GVWl h_{1}'$.
\end{lem}
Deduced from equality $\SIG(h_{1})=\SIG(h_{1}')$ and decreasing $\HM$
during reduction: $\HM(h_{1})\prec_{0}\HM(h_{1}').$
\begin{lem}
\label{lem:greater-or-smaller}Let $h_{1}\GVWl h_{2}$ be labeled
polynomials. Then for $\forall h_{3}\in H\setminus\left\{ \left(0,0\right)\right\} $
at least one of the following two inequalities holds: $h_{1}\GVWl h_{3}$
or $h_{3}\GVWl h_{2}$.
\end{lem}
The lemma clause gives inequality
\begin{equation}
\HM(h_{1})\SIG(h_{2})\prec_{0}\HM(h_{2})\SIG(h_{1})\label{eq:gvw-order-3}
\end{equation}
which shows $\HM(h_{2})\neq0,\SIG(h_{1})\neq0$. Therefore for the
special case $\HM(h_{3})=0$ we get $h_{3}\GVWl h_{2}$and for the
case $\SIG(h_{3})=0$ we get $h_{1}\GVWl h_{3}$. For remaining generic
non-zero case the inequality \eqref{eq:gvw-order-3} can be multiplied
by non-zero monomial $\HM(h_{3})\SIG(h_{3})$: 
\begin{equation}
\HM(h_{3})\SIG(h_{3})\HM(h_{1})\SIG(h_{2})\prec_{0}\HM(h_{3})\SIG(h_{3})\HM(h_{2})\SIG(h_{1}).\label{eq:left-or-right-monom}
\end{equation}
So, the element $\HM(h_{3})^{2}\SIG(h_{2})\SIG(h_{1})\in\mathbb{T}_{0}$
need to be $\succ_{0}$ than left side or $\prec_{0}$ than right
side of inequality \eqref{eq:left-or-right-monom}, and gives after
cancellation one of the inequalities from the lemma statement.

\section{Algorithm}
\begin{lyxlist}{00.00.0000}
\item [{Input:}] polynomial set $\{g_{1},\dots,g_{m}\}$ being a Gröbner
basis; polynomial $f$.
\item [{Variables:}] $R$ and $B$ -- subsets of $H$; $(\sigma,p')\in H$
-- current step's labeled polynomial before reduction; $(\sigma,p)$
-- the same after reduction; $r,b$ -- elements of $R$ and $B$
\item [{Result:}] Gröbner basis of ideal$I=\left(g_{1},\dots,g_{m},f\right)$
\end{lyxlist}

\paragraph*{\label{par:SimpleSignatureGroebner}SimpleSignatureGroebner$\left(\left\{ g_{1},\dots,g_{m}\right\} ,f\right)$}
\begin{enumerate}
\item $R\leftarrow\{\left(\HM(g_{1}),0\right),\left(\HM(g_{2}),0\right),\dots,\left(\HM(g_{m}),0\right),(0,g_{1}),(0,g_{2}),\dots,(0,g_{m})\}$
\item $B\leftarrow\{\}$
\item $(\sigma,p')\leftarrow(1,f)$
\item \textbf{do forever:}

\begin{enumerate}
\item \label{enu:-before-reduce}$p\leftarrow$ReduceCheckingSignatures($\sigma,p',R$)
\item \label{enu:-after-reduce}$R\leftarrow R\cup\left\{ \left(\sigma,p\right)\right\} $
\item \textbf{if} $p\not=0$\textbf{:}

\begin{enumerate}
\item \textbf{for $\{r\in R\,|\, r\GVWl\left(\sigma,p\right),\HM(r)\neq0\}$:}

\begin{enumerate}
\item $B\leftarrow B\cup\{\frac{\LCM(\HM(r),\HM(p))}{\HM(r)}r\}$
\end{enumerate}
\item \textbf{for }$\{r\in R\,|\, r\GVWg\left(\sigma,p\right)\}$\textbf{:}

\begin{enumerate}
\item $B\leftarrow B\cup\{\frac{\LCM(\HM(r),\HM(p))}{\HM(p)}\left(\sigma,p\right)\}$
\end{enumerate}
\end{enumerate}
\item \label{enu:-remove-from-B}$B\leftarrow B\setminus\{b\in B\,|\,\exists r\in R\, r\GVWl b\wedge\SIG(r)|\SIG(b)\}$
\item \textbf{if $B\neq\varnothing$: }$(\sigma,p')\leftarrow$ element
of $B$ with $\prec$-minimal signature
\item \textbf{else: break}
\end{enumerate}
\item \textbf{return} $\{\poly(r)\,|\, r\in R\}$
\end{enumerate}

\paragraph*{\label{par:ReduceCheckingSignatures:}ReduceCheckingSignatures$\left(\sigma,p,R\right)$}
\begin{enumerate}
\item \textbf{do while $\exists r\in R\, r\GVWg(\sigma,p)\wedge\HM(r)|\HM(p)$:}

\begin{enumerate}
\item $p\leftarrow$signature-safe reduce $p$ by $\GVWg$-maximal element
$r$ from the set in cycle clause
\end{enumerate}
\item \textbf{return} $p$\end{enumerate}
\begin{lem}
All pairs from $\mathbb{T}_{0}\times P$ appeared in the algorithm
are labeled polynomials from$H\setminus\left\{ \left(0,0\right)\right\} $.
\end{lem}
The elements created before entering main cycle are labeled polynomials
mentioned above as examples. All other labeled polynomials in the
algorithm are created either with multiplication by $t\in\mathbb{T}$
or with signature-safe reduction, so they satisfy the correctness
property and belongs to $H$.

The clauses of cycles extending $B$ enforces the absence in\textbf{
}$B$ elements with zero signature or zero highest monomial. So, $\sigma$
never can be 0 and the only $R$ elements with zero signatures are
$(0,g_{1}),...,(0,g_{m}).$ Any labeled polynomial added to $R$ can
have zero highest monomial but $R$ does not contain zero polynomial
with zero signature.

\section{Algorithm termination}
\begin{lem}
\label{lem:exist-reductor}At the any moment during the algorithm
execution any labeled polynomial from $B$ can be signature-safe reduced
by some element of $R.$
\end{lem}
Labeled polynomials are added to $B$ in a way ensuring existence
at least one possible signature-safe reductor. The pair $(\sigma,p)\in R$
is such reductor for polynomials added in first \textbf{for} cycle,
and $r\in R$ -- for the polynomials added in the second cycle.
\begin{lem}
\label{lem:r-gvw-small}Before reduction of polynomial $p'$ -- at
the step\ref{enu:-before-reduce} of any algorithm iteration -- the
signatures of elements$\left\{ r\in R\,|\, r\GVWl(\sigma,p')\right\} $
does not divide $\sigma$.
\end{lem}
This holds at the first algorithm iteration because $\sigma=1$ and
$R$ does not contain elements with signatures dividing 1. This holds
during next iterations because the existence such elements in $R$
would lead to removal $\left(\sigma,p'\right)$ from $B$ during previous
iterations at the step \ref{enu:-remove-from-B}.
\begin{lem}
\label{lem:r-gvw-big}After reduction of $p'$ to $p$ -- at the step
\ref{enu:-after-reduce} of any algorithm iteration -- the highest
monomials of elements $\left\{ r\in R\,|\, r\GVWg(\sigma,p)\right\} $
does not divide $\HM(p)$.
\end{lem}
The cycle in the \nameref{par:ReduceCheckingSignatures:} stops only
when it achieves $p$ for which there is no such elements in $R$
.
\begin{lem}
\label{lem:adds-really-new}After reduction of $p'$ to $p$ -- at
the step \ref{enu:-after-reduce} of any algorithm iteration -- no
one from $R$ elements has simultaneously a highest monomial dividing
$\HM(p)$ and a signature dividing $\sigma$.
\end{lem}
The lemma \ref{lem:exist-reductor} ensures that $p'$is reduced at
least once, so $(\sigma,p')\GVWg(\sigma,p)$. Now, by lemma\ref{lem:greater-or-smaller}
for $\forall r\in R$ we have $r\GVWg(\sigma,p)$ or $r\GVWl(\sigma,p')$.
These inequalities allow to apply either lemma \ref{lem:r-gvw-small}
or lemma \ref{lem:r-gvw-big}.
\begin{thm}
The algorithm\nameref{par:SimpleSignatureGroebner} terminates
\end{thm}
To prove the termination we need to show that all \textbf{do} cycles
stops after finite number of executions. In the cycle inside \nameref{par:ReduceCheckingSignatures:}
with non-zero $p$ during every iteration we have $\HM(p)$ decrease
according to $\prec_{0}$, which is possible only finite number of
times. When $p$ becomes zero it stops immediately because of $\GVWl$-minimality
of $(\sigma,0)$.

The set $R\subset\mathbb{T}_{0}\times P$ is extended every step of
the main algorithm cycle. It can be splitted to $R_{*0}\cup R_{0*}\cup R_{**},$
where $R_{*0}\subset\mathbb{T}\times\left\{ 0\right\} ,R_{0*}\subset\left\{ 0\right\} \times P\setminus\left\{ 0\right\} ,R_{**}\subset\mathbb{T}\times P\setminus\left\{ 0\right\} $.
$R_{0*}$ does never extend because $\sigma\neq0$. For sets $R_{*0}$
and $R_{**}$ we apply a method based on idea of monoid ideal introduced
in \cite{KreuzerRobbianoBook1} as ``monoideal''. Consider the following
two sets which are monoideals: $L_{*0}=\left(\left\{ \sigma\,|\,(\sigma,0)\in R_{*0}\right\} \right)\subset\mathbb{T}$
and $L_{**}=\left(\left\{ (\sigma,t)\,|\,\exists(\sigma,p)\in R_{**}\, t=\HM(p)\right\} \right)\subset\mathbb{T}\times\mathbb{T}$.
Lemma \ref{lem:adds-really-new} shows that elements being added to
$R$ expand either $L_{*0}$ or $L_{**}$ in every cycle iteration.
The monoids $\mathbb{T}$ and $\mathbb{T}\times\mathbb{T}$ are isomorphic
to $\mathbb{N}^{n}$ and $\mathbb{N}^{2n}$, so the Dickson's lemma
can be applied to their monoideals. It states exactly the needed fact
-- only finite number of expansions is possible for such monoideals.

\section{Correctness of output}
\begin{defn}
\emph{S-representation of} $h\in H$ over set $\left\{ r_{i}\right\} \subset H$
is an expression $\poly(h)=\sum_{j}K_{j}t_{j}\poly(r_{i_{j}}),$ $K_{j}\in k,t_{j}\in\mathbb{T},i_{j}\in\mathbb{N}$,
such that $\forall j\,\HM(h)\succcurlyeq_{0}\HM(t_{j}r_{i_{j}}),\SIG(h)\succcurlyeq_{0}\SIG(t_{j}r_{i_{j}})$.\end{defn}
\begin{lem}
\label{lem:one-j-exact}Let $\poly(h)=\sum_{j}K_{j}t_{j}\poly(r_{i_{j}})$
be S-representation of $h$. Then at least one $j$ satisfies $\HM(h)=\HM(t_{j}r_{i_{j}})$.
\end{lem}
To get a $j$ satisfying the equality we can take a value which gives
the $\succ$-maximum of $\HM(t_{j}r_{i_{j}})$.

The next definition extends the notation of S-basis from \cite{TheF5Revised}:
\begin{defn}
We call a labeled polynomial set $R\subset H$\emph{ S-basis} (correspondingly
\emph{S$_{\sigma}$-basis}), if all elements of $H$ (correspondingly
$\left\{ h\in H\,|\,\SIG(h)\prec_{0}\sigma\right\} $)\emph{ }have
S-representation over $R$\emph{.}\end{defn}
\begin{lem}
\label{lem:s-basis-and-no-reductions}Let $\sigma\succ_{0}0,R=\left\{ r_{i}\right\} $
be S\emph{$_{\sigma}$}-basis and $h_{1},h_{2}\in H,\SIG(h_{i})=\sigma$
be labeled polynomials, that can't be signature-safe reduced by $R$
elements. Then $\HM(h_{1})=\HM(h_{2})$ and $h_{1}$ has an S-representation
over $R\cup\left\{ h_{2}\right\} $.
\end{lem}
We have from the definition of $H$ that $\exists u_{i}\in P\,\HM(u_{i})=\sigma,u_{i}f\equiv\poly(h_{i})\pmod{I_{0}},i=1,2.$
It means that there exists a linear combination of $\poly(h_{i})$
having signature $\prec_{0}\sigma$. This can be written as: 
\[
\exists K\in k,v\in P\,\HM(v)=\sigma'\prec_{0}\sigma,vf\equiv\poly(h_{1})-K\poly(h_{2})\pmod{I_{0}},
\]
or in the terminology of labeled polynomials: $\left(\sigma',p'\right)=\left(\sigma',\poly(h_{1})-K\poly(h_{2})\right)\in H$.
From the definition of S\emph{$_{\sigma}$}-basis and the property
$\sigma'\prec_{0}\sigma$ we conclude: $\exists r_{j}\in R,t\in\mathbb{T}\,\SIG(tr_{j})\preccurlyeq_{0}\sigma',\HM(tr_{j})=\HM(p')$.
So $\HM(h_{i})\neq\HM(p'),i=1,2$, because in the case of equality
$r_{j}$ would be signature-safe reductor for $h_{i}$. It is possible
only if $\HM(h_{i})$ are canceled while subtraction with $k$-coefficient,
what means that $\HM(h_{1})=\HM(h_{2})$. S-representation of $h_{1}$
is constructed by adding $K\poly(h_{2})$ to S-representation of $\left(\sigma',p'\right)$.
\begin{thm}
\label{thm:exist-r-sigma}Every iteration of the algorithm after step
\ref{enu:-remove-from-B} the following invariant holds: for $\forall\sigma\in\mathbb{T},\sigma\prec$
signatures of elements of $B$, exists $r_{\sigma}\in R,t_{\sigma}\in\mathbb{T}:\SIG(t_{\sigma}r_{\sigma})=\sigma$
such that $t_{\sigma}r_{\sigma}$ can't be signature-safe reduced
by $R$.
\end{thm}
The set $R_{\sigma}=\left\{ r\in R\,|\,\SIG(r)|\sigma\right\} $ is
not empty, because contains the element $r_{0}$ added during the
first algorithm iteration with $\SIG(r_{0})=1$. Let $r_{\sigma}$
be $\GVWl$-minimal element of the set; take $t_{\sigma}=\frac{\sigma}{\SIG(r_{\sigma})}$.
Suppose that $t_{\sigma}r_{\sigma}$ can be signature-safe reduced
by some $r_{1}\in R$. This gives that $r_{1}\GVWg r_{\sigma}$ and
both sides of inequality are non-zero. It means that during the iteration
which inserts in $R$ the last of $\{r_{\sigma},r_{1}\}$ the set
$B$ was extended by labeled polynomial $t'r_{\sigma}$, where $t'=\frac{\LCM(\HM(r_{1}),\HM(r_{\sigma}))}{\HM(r_{\sigma})}$
and $t'|t_{\sigma}$. So we have $\SIG(t'r_{\sigma})|\SIG(t_{\sigma}r_{\sigma})=\sigma\Rightarrow\SIG(t'r_{\sigma})\preccurlyeq\sigma\prec$
signatures of elements of $B$. This signatures inequality implies
that $t'r_{\sigma}$ can't be element of $B$ during the current iteration
and was removed at the step \ref{enu:-remove-from-B} of some previous
iteration, so $\exists r_{2}\in R\, r_{2}\GVWl t'r_{\sigma},\SIG(r_{2})|\SIG(t'r_{\sigma})$.
This is impossible, because the existence of $r_{2}\GVWl r_{\sigma},r_{2}\in R_{\sigma}$
contradicts $\GVWl$-minimality of $r_{\sigma}$.
\begin{thm}
\label{thm:has-s-repr}Every iteration of the algorithm after step
\ref{enu:-remove-from-B} the following invariant holds: $\forall h\in H,\SIG(h)\prec$
signatures of elements of $B$ has S-representation over $R$.
\end{thm}
Suppose that invariant breaks during some algorithm iteration and
take the $\prec_{0}$-minimal $\sigma$ that has non-empty corresponding
set $V_{\sigma}\eqdef\{h\in H\,|\, h\mbox{ breaks invariant},\SIG(h)=\sigma\}$.
Then $R$ is S\emph{$_{\sigma}$}-basis. $\forall g\in I_{0}\,\left(0,g\right)$
has S-representation over $\{(0,g_{1}),...,(0,g_{m})\}\subset R$,
so $\sigma\succ_{0}0$. Select $v_{\sigma}$ -- one of the $V_{\sigma}$
elements with $\prec_{0}$-minimal $\HM$. It can't be signature-safe
reduced by $R$ because the reduction result $v_{1}$ would be element
of$V_{\sigma}$ with $\HM(v_{1})\prec_{0}\HM(v_{\sigma})$. Take $w_{\sigma}\eqdef t_{\sigma}r_{\sigma}$
from the invariant of theorem \ref{thm:exist-r-sigma} and apply lemma
\ref{lem:s-basis-and-no-reductions} to $v_{\sigma},w_{\sigma}$ and
$R$. The lemma says that $v_{\sigma}$ has S-representation over
$R\cup\{w_{\sigma}\}$. All entries of $w_{\sigma}$ in the representation
can be replaced by $t_{\sigma}r_{\sigma}$ to acquire S-representation
of$v_{\sigma}$ over $R$ only. It's existence leads to contradiction.
\begin{lem}
If $R$ is S-basis, then $\{\poly(r)\,|\, r\in R\}$ is a Gröbner
basis of ideal $I$.
\end{lem}
For $\forall p\in I$ we can take some $h=\left(\sigma,p\right)\in H$
and apply lemma \ref{lem:one-j-exact} to it.
\begin{thm}
\nameref{par:SimpleSignatureGroebner} returns Gröbner basis
\end{thm}
At the moment of algorithm termination $B=\varnothing$ so by theorem
\ref{thm:has-s-repr} $R$ is S-basis.

\section{Comparison with other algorithms}

The presented algorithm belongs to the family of the Gröbner basis
algorithms using signatures and being at some degree a modification
of F5 algorithm from \cite{FaugereF5}. One of the main modification
directions of F5 is simplifying and clarifying the connected theory
usually bound with some thoughts about extending the area of inputs
the algorithm can be applied to. Investigations in this direction
can be found in \cite{GermanF5Proof,ZobninGeneralization,F5InBBStyle}.
The other direction is improving efficiency of computations by introducing
criteria to detect and don't perform some unnecessary computations.
It is studied in \cite{F5C,G2V,SignatureBasedGBs} and allows to perform
computations in a way that reduces to the end only polynomials that
are either new S-basis entries or corresponds to extending monoideal
of signatures known to be zero polynomial signatures -- called \emph{syzygy
signatures.} Generalization simultaneously using all criteria described
in algorithms TRB-MJ and SB \cite{HuangConception,PracticalGB} achieve
even more efficiency because all discardings are performed before
any computational heavy operations like polynomial reduction or computation
highest monomial of S-polynomial -- so the non-trivial computations
are never become unnecessary because their results are never discarded.

All mentioned algorithms including original F5 use discarding criteria
of two types: syzygy-based criteria and rewrite-like criteria, with
separate proof of correctness for each type. The other common idea
used in the algorithms are S-polynomials: even the algorithms that
does not deal with S-polynomials directly make heavy usage of them
in the correctness proof.

This paper describes an algorithm computing minimal S-basis and discarding
computations with efficiency identical to TRB-MJ, but using the only
one discarding criteria in step \ref{enu:-remove-from-B}, which is
based on $\GVWl$-ordering of $R$. The routine \nameref{par:ReduceCheckingSignatures:}
can use different reductor selection strategies and their efficiency
is open question. The method proposed in this paper is based on the
same ordering of $R$ and is identical for homogeneous case with the
methods used in original F5. The correctness proof is given without
the use of S-polynomials and is formulated in a way that allows to
apply the clear algebraic interpretation of signature-based algorithms
from \cite{ZobninGeneralization} to the presented algorithm.

Algorithm simplification lead to simplification of programming its
implementation and debugging. It is achieved by the smaller number
of objects involved in computation and use of the same order for discarding
criteria and reductor selection. Simplicity of implementation and
the absence of complex data structures allows quick algorithm integration
with any computer algebra system that can work with polynomials. The
author's implementation linked below was written from scratch in a
8 hours what is a lot smaller than the time author spent implementing
other algorithms with a similar tools.

The algorithm was implemented in C++ using low-level functions from
computer algebra system Singular 3-1-4 and open source codes of C.
Eder (one of the authors of \cite{SignatureBasedGBs}) for implementing
F5-like algorithms in this system. The source is contained in ``ssg''
function in a file available at \href{https://github.com/galkinvv/Singular-f5-like/blob/ssg/kernel/kstd2.cc}{https://github.com/galkinvv/Singular-f5-like/blob/ssg/kernel/kstd2.cc}

Comparison of SimpleSignatureGroebner implementation with other Gröbner
basis algorithms implemented by C. Eder gives practical checks for
the following theoretical facts:
\begin{itemize}
\item algorithm SimpleSignatureGroebner correctly computes Gröbner basis;
\item the number of polynomials in the result set is not greater than the
number of polynomials in a result of other incremental algorithms
that compute S-basis;
\item the execution time is not greater than execution time of other signature-based
incremental algorithms.
\end{itemize}
\novestnikonly{

\bibliographystyle{plain}
\bibliography{f5_references}

\begin{thebibliography}{10}

\bibitem{TheF5Revised}
A.~{Arri} and J.~{Perry}.
\newblock The f5 criterion revised.
\newblock {\em ArXiv e-prints}, December 2010.

\bibitem{F5C}
C.~{Eder} and J.~{Perry}.
\newblock F5c: a variant of faug\`{e}re's f5 algorithm with reduced gr\"{o}bner
  bases.
\newblock {\em ArXiv e-prints}, June 2009.

\bibitem{SignatureBasedGBs}
C.~{Eder} and J.~{Perry}.
\newblock Signature-based algorithms to compute gr\"{o}bner bases.
\newblock {\em ArXiv e-prints}, January 2011.

\bibitem{FaugereF5}
Jean~Charles Faug\`{e}re.
\newblock A new efficient algorithm for computing gr\"{o}bner bases without
  reduction to zero (f5).
\newblock In {\em Proceedings of the 2002 international symposium on Symbolic
  and algebraic computation}, ISSAC '02, pages 75--83, New York, NY, USA, 2002.
  ACM.

\bibitem{G2V}
Shuhong Gao, Yinhua Guan, and Frank Volny, IV.
\newblock A new incremental algorithm for computing gr\"{o}bner bases.
\newblock In {\em Proceedings of the 2010 International Symposium on Symbolic
  and Algebraic Computation}, ISSAC '10, pages 13--19, New York, NY, USA, 2010.
  ACM.

\bibitem{GermanF5Proof}
O.~{German}.
\newblock Proof of the faug\`{e}re criterion for the f5 algorithm.
\newblock {\em Mathematical Notes}, 88(4):502--510, 2010.

\bibitem{HuangConception}
L.~{Huang}.
\newblock A new conception for computing gr\"{o}bner basis and its
  applications.
\newblock {\em ArXiv e-prints}, December 2010.

\bibitem{KreuzerRobbianoBook1}
Martin {Kreuzer} and Lorenzo {Robbiano}.
\newblock {\em Computational commutative algebra. 1}.
\newblock Springer-Verlag, Berlin, 2000.

\bibitem{PracticalGB}
B.~{Roune} and M.~{Stillman}.
\newblock Practical gr\"{o}bner basis computation.
\newblock 2012.

\bibitem{F5InBBStyle}
Y.~{Sun} and D.~{Wang}.
\newblock The f5 algorithm in buchberger's style.
\newblock {\em ArXiv e-prints}, June 2010.

\bibitem{ZobninGeneralization}
A.~I. Zobnin.
\newblock Generalization of the f5 algorithm for calculating gr\"{o}bner bases
  for polynomial ideals.
\newblock {\em Program. Comput. Softw.}, 36(2):75--82, March 2010.

\end{thebibliography}

}
\end{document}